\documentclass[12pt]{amsart}
\usepackage{amsmath, amsthm, amssymb}
\usepackage{amssymb}
\usepackage{amsfonts}
\usepackage{amscd}
\usepackage{mathrsfs}
\usepackage{latexsym}
\usepackage{amsmath}
\usepackage{txfonts}
\usepackage{tikz}
\usepackage{graphicx}
\usepackage{subfig}
\usepackage{tikz-cd}
\usepackage{blindtext}

\usepackage[all, cmtip]{xy}

\newtheorem{theorem}{Theorem}[section]
\newtheorem{lemma}[theorem]{Lemma}
\newtheorem{corollary}[theorem]{Corollary}
\newtheorem{proposition}[theorem]{Proposition}

\theoremstyle{definition}
\newtheorem{definition}[theorem]{Definition}
\newtheorem{example}[theorem]{Example}

\theoremstyle{remark}
\newtheorem{claim}[theorem]{Claim}

\newcommand{\spec}[1]{\textrm{spec}#1}
\newcommand{\sshf}[1]{\mathscr{O}_{#1}}
\newcommand{\shf}[1]{\mathscr{#1}}

\newcommand{\iso}{\simeq}
\newcommand{\ses}[3]{0\rightarrow#1\rightarrow#2\rightarrow#3\rightarrow{0}}

\newcommand{\is}[1]{\mathscr{I}_{#1}}

\newcommand{\paren}[1]{\left(#1\right)}

\numberwithin{equation}{section}

\begin{document}
\allowdisplaybreaks
\title[On direct images of twisted pluricanonical sheaves]{On direct images of twisted pluricanonical sheaves on normal varieties}

\author{Chih-Chi Chou}
\author{Lei Song\**}

\address{Seattle, WA 98107, United States}
\email{cchou9@gmail.com}
\address{School of Mathematics, Sun Yat-sen University, Guangzhou, Guangdong 521075, China}
\email{songlei3@mail.sysu.edu.cn}


\dedicatory{}


\begin{abstract}
We study the depth properties of certain direct image sheaves on normal varieties. Let $f: Y\rightarrow X$ be a proper morphism of relative dimension $d$ from a smooth variety onto a normal variety such that the preimage $E$ of the singular locus of $X$ is a divisor. We show that for any integer $m>0$, the higher direct image $R^df_*\omega^{\otimes m}_Y(aE)$ modulo the torsion subsheaf is $S_2$, provided that $a$ is sufficiently large. In case $f$ is birational, we give criteria on $a$ for the direct image $f_*\omega_Y(aE)$ to coincide with $\omega_X$. We also introduce an index measuring the singularities of normal varieties.
\end{abstract}

\maketitle

\section{Introduction}
Let $X$ be a normal proper variety over a field of characteristic zero and $\mu: Y\rightarrow X$ be a log resolution of singularities such that $\mu$ is isomorphic over the smooth locus of $X$. Let $E$ be the full exceptional divisor of $\mu$, which is reduced. There exists a non-decreasing sequence of coherent sheaves (cf. \S 2.2)
\begin{equation}\label{increasing sequence}
    0\subseteq \mu_*\omega_{Y}\subseteq \mu_*\omega_{Y}(E)\subseteq \mu_*\omega_{Y}(2E)\subseteq\cdots \mu_*\omega_{Y}(aE)\subseteq\cdots \omega_X,
\end{equation}
where $\omega_X$ denotes the dualizing sheaf. By the coherence of $\omega_X$, one can show that  (\ref{increasing sequence}) stabilizes at $\omega_X$ for arbitrary singularities, see e.g.~Proposition \ref{Main}. 

When $X$ has at worst rational singularities, $\mu_*\omega_Y\iso \omega _X$. And when $X$ has at worst Du Bois singularities, $\mu_*\omega_Y(E)\iso \omega_X$. If in addition $X$ is Cohen-Macaulay, then the above isomorphisms indeed characterize rational singularities by G. Kempf , and Du Bois singularities by Kov\'{a}cs-Schwede-Smith \cite{KSS} respectively. 

The purpose of the article is twofold. On the one hand, we introduce in Definition \ref{index def} an index of the singularities of $X$ measuring how quickly the sequence (\ref{increasing sequence}) stabilizes to $\omega_X$ among all log resolutions. A necessary condition (Proposition \ref{necessary condition}) and a sufficient condition (Theorem \ref{splitting theorem}) are given for when the direct images coincide with $\omega_X$. Suppose the log resolution $\mu$ resolves the Jacobian ideal of $X$, as a consequence of Theorem \ref{direct image}, we show in Corollary \ref{MJ Discrepancy} that one can tell when (\ref{increasing sequence}) stabilizes in terms of Mather-Jacobian discrepancies. All these results may be viewed as giving certain upper or lower bound for the index of singularities.

On the other hand, an alternative way to interpret that the direct images in (\ref{increasing sequence}) coincides with $\omega_X$ for large enough number $a$ is that these sheaves satisfy Serre's $S_2$ condition, and hence are reflexive. From this point of view, we can generalize Proposition \ref{Main} to arbitrary relative dimension (Theorem \ref{Main Theorem}). We first prove the case $m=1$ by induction on relative dimension, then obtain the result for pluri-forms by Viehweg's cyclic covering trick.

Finally it is worth noting that we do not assume $X$ is $\mathbb{Q}$-Gorenstein throughout the article.

Conventions: We work over an algebraically closed field of characteristic zero. By a variety, we mean a separated integral scheme of finite type over the base field. For a given coherent sheaf $\shf{F}$, $\shf{F}/{\text{tor}}$ always denotes the quotient of $\shf{F}$ by its torsion subsheaf.

{\it Acknowledgments.}
We would like to thank Lawrence Ein, Hailong Dao, S\'andors Kov\'acs, Bangere Purnaprajna and Qi Zhang for the helpful discussions and suggestions during the preparation of the paper. We also thank Linquan Ma for answering our questions. Finally, we are grateful to the referee for careful reading of the article and valuable suggestions.

\section{Preliminaries}

\subsection{$S_k$ conditions and reflexive sheaves}(cf. \cite{Hartshorne80})
Let $X$ be a Noetherian scheme. Given a coherent sheaf $\shf{F}$ on $X$, we say $\shf{F}$ satisfies Serre's $S_k$ condition if $\text{depth}_{\frak{m}_x}{F_x}\ge \min\{k, \dim \sshf{X, x}\}$ for any $x\in X$, where $(\sshf{X, x}, \frak{m}_x)$ is the local ring of $x$. We will only need $S_2$ condition in this article. We denote the dual sheaf $\mathscr{H}om_{X}(\shf{F}, \sshf{X})$ of $\shf{F}$ by ${\shf{F}}^{\vee}$. There is a functorial natural morphism of coherent sheaves $\shf{F}\rightarrow \shf{F}^{\vee\vee}$. For any integer $m\ge 1$, $(\shf{F}^{\otimes m})^{\vee\vee}$ is denoted by $\shf{F}^{[m]}$. A sheaf $\shf{F}$ is called \textit{reflexive} if the morphism is an isomorphism. When $X$ is normal, a torsion free coherent sheaf $\shf{F}$ on $X$ is reflexive if and only if it is $S_2$. In particular, $\omega_X$ is reflexive.

We collect a few useful lemmas on depth that are needed in the sequel.

\begin{lemma}\label{S_k condition}
Let $X$ be a Noetherian scheme. For any coherent sheaves $\shf{F}$, $\shf{G}$ and invertible sheaf $\shf{L}$. It hold that
\begin{itemize}
  \item $\shf{F}\oplus\shf{G}$ is $S_k$ $\iff$ both $\shf{F}$ and $\shf{G}$ are $S_k$.
  \item $\shf{F}$ is $S_k$ $\iff $ $\shf{F}\otimes_{\sshf{X}}\shf{L}$ is $S_k$.
\end{itemize}
\end{lemma}
\begin{proof}
These follow from local cohomology criterion for depth, see \cite[Cor.~3.10]{Hartshorne67} and that $\shf{L}_x\iso \sshf{X, x}$ for all $x\in X$.
\end{proof}

\begin{lemma}\label{criterion1}
Let $\ses{\shf{F}}{\shf{G}}{\shf{Q}}$ be a short exact sequence of coherent sheaves on a normal scheme $X$. Suppose $\shf{G}$ is reflexive and the support of $\shf{Q}$ has codimension $\ge 2$. Then $\shf{F}\iso \shf{G}$ if and only if $\shf{F}$ is reflexive.
\end{lemma}
\begin{proof}
This is a direct consequence of \cite[Cor.~1.5]{Hartshorne80}.
\end{proof}

\begin{lemma}(\cite[Prop.~5.4]{KollarMori08})\label{depth under finite map}
Let $f: Y\rightarrow X$ be a finite morphism of Noetherian schemes and $\shf{F}$ a coherent sheaf on $Y$ such that $\text{Supp}(\shf{F})$ is pure dimensional. Then $\shf{F}$ is $S_k$ if and only if $f_*\shf{F}$ is $S_k$. \qed
\end{lemma}

\subsection{Filtrations of $\omega^{[m]}_X$ associated to maps}

Let $X$ be a proper normal variety and $\omega_X$ the dualizing sheaf. The singular locus $X_{\text{sing}}$ of $X$ has codimension $\ge 2$. Denote by $U=X\backslash X_{\text{sing}}$ the smooth locus, and by $j: U\rightarrow X$ the open immersion. Let $f: Y\rightarrow X$ be a proper birational morphism from a smooth variety Y such that $f|_V: V=f^{-1}(U)\rightarrow U$ is an isomorphism.

Fix an integer $m\ge 1$. For any integer $a\ge 0$ and an integral divisor $E$ supported on the exceptional locus of $f$ ($\text{supp}(E)$ is possibly smaller than $\text{Exc}(f)$), there are natural morphisms
\begin{equation*}
    f_*\omega^{\otimes m}_{Y}(aE)\rightarrow j_*j^*f_*\omega^{\otimes m}_{Y}(aE)\xrightarrow{\iso}j_*(f|_V)_*\paren{\omega^{\otimes m}_{Y}(aE)|_V}\iso j_*\omega^{\otimes m}_U\xleftarrow{\iso}\omega^{[m]}_X,
\end{equation*}
where the rightmost map is an isomorphism because $\omega^{[m]}_X$ is reflexive \cite[Section 1]{Hartshorne80}. The above composite gives a natural morphism of coherent sheaves
\begin{equation*}
    f_*\omega^{\otimes m}_{Y}(aE)\rightarrow \omega^{[m]}_X.
\end{equation*}
Since the sheaf $f_*\omega^{\otimes m}_{Y}(aE)$ is torsion free of rank 1 and $f_*\omega^{\otimes m}_{Y}(aE)\rightarrow \omega^{[m]}_X$ is an isomorphism away from $X_{\text{sing}}$, the natural map $f_*\omega^{\otimes m}_{Y}(aE)\hookrightarrow \omega^{[m]}_X$ is injective. Therefore one has the increasing sequence
\begin{equation*}
    0\subseteq f_*\omega^{\otimes m}_{Y}\subseteq f_*\omega^{\otimes m}_{Y}(E)\subseteq f_*\omega^{\otimes m}_{Y}(2E)\subseteq\cdots \omega^{[m]}_X,
\end{equation*}
which stablizes at $\omega^{[m]}_X$ by the following

\begin{proposition}\label{Main}
Let $f:Y \rightarrow X$ be a proper birational morphism between Noetherian normal schemes. Let $E$ denote the reduced full exceptional divisor and assume $E$ is Cartier. Then for any integer $m\ge 1$, we have
$f_*\omega^{[m]} _Y(aE)\iso \omega^{[m]} _X$ for $a\gg 0$.\qed
\end{proposition}
The proof is based on the Noetherian property and the fact that any section of $\omega^{[m]} _X$  has a finite order of pole along any exceptional prime divisor. As it is well-known to the experts, we shall omit it here.

In view of  Proposition \ref{Main} and Lemma \ref{criterion1}, $f_*\omega^{\otimes m}_Y(aE)$ is reflexive when $a\gg 0$.  Bearing this in mind, we consider the slightly more general situation.

Let $f: Y\rightarrow X$ be a surjective proper morphism from a smooth variety $Y$ to a normal variety $X$ of relative dimension $d=\dim Y-\dim X$. Suppose $E=f^{-1}(X_{\text{sing}})$ is a divisor, with the reduced induced structure. As before let $U$ be the smooth locus of $X$ and $V=f^{-1}(U)$. Moreover suppose that $f$ has connected fibres. By \cite[Prop. 7.6]{Kollar86}, \footnote{Though it states for the projective case, one can extend the statement to general case by compactification.} it holds that
\begin{equation*}
    R^df_*\omega_V\iso \omega_U.
\end{equation*}
 For $a\ge 0$, we have a natural map
\begin{equation*}
    R^df_*(\omega_Y(aE))\rightarrow \omega_X.
\end{equation*}
In general when $a>0$, the direct images of twisted pluricanonical sheaves $R^df_*(\omega_Y(aE))$ have the torsion parts. Modulo the torsion parts, one gets the increasing sequence
\begin{equation}\label{filtration}
     0\subseteq R^df_*\omega_Y\subseteq R^df_*\paren{\omega_Y(E)}/{\text{tor}}\subseteq R^df_*\paren{\omega_Y(2E)}/{\text{tor}}\subseteq\cdots \omega_X.
\end{equation}
We shall show in \S 4 that it stablizes at $\omega_X$ ultimately.

\section{Index of Singularities}
In view of Proposition \ref{Main}, it is natural to introduce the following index measuring the singularities of a normal variety.
\begin{definition}\label{index def}
Given a normal variety $X$, the index of $X$ is defined as
\begin{equation*}
m_X=\min _{f:Y\rightarrow X}\{a\},
\end{equation*}
where the minimum is taken over all log resolution of $X$ with the property that $f_*\omega_Y(aE_Y)\iso \omega_X$.
\end{definition}

\example(Generalized affine cones)
Let $X$ be a smooth projective variety of dimension $d$ and $L$ an ample line bundle on $X$ such that $K_X$ and $L$ are linearly independent in $NS(X)_{\mathbb{Q}}$. Consider the section ring
 \begin{equation*}
    R(X, L)=\oplus_{m\ge 0} H^0(X, L^m).
 \end{equation*}
Then the affine cone $V=\spec R(X, L)$ is a normal algebraic variety. In case $L$ is very ample, it is the normalization of the classical affine cone associated to the embedding by $|L|$. By our assumption, $V$ is not $\mathbb{Q}$-Gorenstein, cf. \cite[Prop. 3.14]{Kollar13}.

We have a birational morphism
\begin{equation*}
    p: V':=\spec_{\sshf{X}}\paren{\oplus_{m\ge 0} L^m}\rightarrow V,
\end{equation*}
with the exceptional divisor $E\iso X$. Since $V'$ is the total space of $L^{-1}$, $V'$ is smooth.

For any integer $a\ge 0$, $p_*\omega_{V'}(aE)$ differs from $\omega_V$ only at the vertex $0$, which corresponds to the irrelevant maximal ideal of $R(X, L)$.
Note that $\sshf{E}(-E)\iso \sshf{E}(1)$ is identified with $L$ via the zero section of $\pi: V'\rightarrow X$. Using the theorem of formal functions (cf. \cite[III Thm 11.1]{Hartshorne77}), we have
\begin{equation*}
   \widehat{p_*\omega_{V'}(aE)_0}\iso \oplus_{m>0} H^0(X, \omega_X\otimes L^{m-a}).
\end{equation*}
On the other hand, by a slight abuse of notation and using \cite{HashimotoKurano11}, we have
\begin{equation*}
    \widehat{{\omega_V}}_0=\oplus_{m\in\mathbb{Z}}H^0(X, \omega_X\otimes L^m).
\end{equation*}
Therefore we see that
\begin{equation*}
    p_*\omega_{V'}(aE)\iso \omega_V, \quad\quad \text{if } a\ge m_V,
\end{equation*}
where
\begin{eqnarray*}
 m_V&=&\min\{ \:a\in\mathbb{Z}_{+} \:|\: H^0(X, \omega_X\otimes L^{-i})=0,  \textit{ for all } i\ge a\} \\
   &=& \min\{ \:a\in\mathbb{Z}_{+} \:|\: H^d(X, L^i)=0,  \textit{ for all } i\ge a\}. 
\end{eqnarray*}

\example
(Secant varieties). Given a smooth projective variety $X$ of dimension $d$ and a sufficiently positive adjoint line bundle $L$ on $X$. For the secant variety $\Sigma=\Sigma(X, L)$, there is a natural resolution of singularities $f: Y\rightarrow \Sigma$ with the exceptional divisor $E$. By \cite{ChouSong17}, $\Sigma$ has Du Bois singularities and it always holds that $f_*\omega_Y(E)\iso \omega_{\Sigma}$, so $m_{\Sigma}\le 1$. Moreover $f_*\omega_Y\iso \omega_{\Sigma}$ if and only if $H^d(X, \sshf{X})=0$.

In general, the index in Definition \ref{index def} is hard to compute, as it is presumably computed over all log resolutions. The following result may be viewed as to give a lower bound on such index.

\begin{proposition}\label{necessary condition}
Let $W=X_{\text{sing}}$ be with the reduced induced scheme structure. Suppose for a log resolution $f: Y\rightarrow X$ and $m\ge 1$, $f_*\omega^{\otimes m}_Y(aE)\iso \omega^{[m]}_X$ holds for some $a$. Then
\begin{equation*}
\is{W}^{a} \cdot \omega^{[m]}_X \subseteq (f_*\omega_Y)^{[m]}.
\end{equation*}
\end{proposition}

\begin{proof}
We give the proof for $m=1$, as that for general $m$ is similar. To begin with, we have the sequence
\begin{equation*}
f_*\omega_Y \subseteq f_*\omega_Y(E)\subseteq \cdots \subseteq f_*\omega_Y((a-1)E)
\subseteq f_*\omega_Y(aE)=\omega_X,
\end{equation*}
and to prove the inclusion statement, it suffices to show that
\begin{equation}\label{injectivity}
\is{W}\cdot f_*\omega_Y(kE)\hookrightarrow f_*\omega_Y((k-1)E), \quad \text{ for } 1\le k  \le a.
\end{equation}
To this end, note that $\is{W}\cdot f_*\omega_Y(kE)$ is torsion free, because
it is the image of the natural map
\begin{equation*}
\is{W}\otimes f_*\omega_Y(kE)\rightarrow f_*\omega_Y(kE),
\end{equation*}
and hence a subsheaf of the torsion free sheaf $f_*\omega_Y(kE)$. So we just need to show that there exists a natural nonzero map $\is{W}\cdot f_*\omega _Y(kE)\rightarrow f_*\omega_Y((k-1)E)$, then the injectivity of (\ref{injectivity}) follows.

We claim that there is a natural map $\rho: f_*\omega _Y(kE) \big| _W \rightarrow f_* \left(\omega_Y(kE)\big| _E\right)$. In fact, let $Y_W$ be the fibre of $W$. By the setup, $E$ is the reduced scheme associated to $Y_W$; so there is
a closed embedding $E \hookrightarrow Y_W$, inducing a map $ \omega_Y(kE)\big| _{Y_W} \rightarrow \omega_Y(kE)\big| _E$.
Pushing it down to $X$, we get the following composition,
\begin{equation*}
\rho: f_*\omega _Y(kE) \big| _W \rightarrow  f_*\left(\omega_Y(kE)\big| _{Y_W}\right) \rightarrow f_*\left(\omega_Y(kE)\big| _E\right),
\end{equation*}
where the first arrow is natural by fibre product, cf.~\cite[III Remark 9.3.1]{Hartshorne77}.

Thus we arrive at the following commutative diagram with exact rows

$$\xymatrix{
 0\ar[r] & \is{W}\cdot f_*\omega_Y(kE)\ar[r]^{\alpha}\ar@{-->}[d]^{\exists} & f_*\omega _Y(kE) \ar[r]\ar@{=}[d] & f_*\omega _Y(kE) \big| _W\ar[r]\ar[d]^{\rho} & 0 \\
 0\ar[r] & f_*\omega_Y((k-1)E)\ar[r] &  f_*\omega _Y(kE)\ar[r]^{\beta} & f_* (\omega_Y(kE)\big| _E). & }. $$
Therefore $im(\alpha)$ is contained in $ker(\beta)$, and hence the first vertical map exists. This finishes the proof.
\end{proof}

In the rest of the section, we focus on the direct images of twisted canonical sheaves. The following theorem gives a sufficient condition for when the direct image coincides with the dualizing sheaf, and hence gives an upper bound of the index of $X$.
\begin{theorem}\label{splitting theorem}
Let $f:Y\rightarrow X$ be a log resolution of a normal variety $X$. Let $E$ be an effective divisor on $Y$ whose
support is contained in the exceptional locus.
If for some integer $a\ge 0$, the natural map $f_* \shf{O}_Y(-aE)\rightarrow \mathcal{R}f_*\shf{O}_Y(-aE)$ in the derived category has a left inverse, then we have
\begin{equation*}
f_*\omega _Y(aE)\iso \omega _X.
\end{equation*}
\end{theorem}

\begin{proof}
Put $\shf{I}= f_*\shf{O}_Y(-aE)$. By assumption we have
\begin{equation*}
\shf{I}\rightarrow \mathcal{R}f_*\shf{O}_Y(-aE) \rightarrow \shf{I}
\end{equation*}
and the composition is $\text{id}_{\shf{I}}$. Applying the functor $\mathcal{R}Hom_{X}( \:, \omega _{X}^{\bullet})$, we get
\begin{equation*}
\mathcal{R}Hom_X(\shf{I}, \omega _{X}^{\bullet})\rightarrow
\mathcal{R}Hom_X(\mathcal{R}f_*\shf{O}_Y(-aE), \omega _{X}^{\bullet})\rightarrow
\mathcal{R}Hom_X(\shf{I}, \omega _{X}^{\bullet})
\end{equation*}
 Applying Grothendieck duality to the middle term and taking the $-n$ cohomology yield
\begin{equation}\label{split of ext}
\mathcal{E}xt^{-n}(\shf{I}, \omega_{X}^{\bullet})
\rightarrow f_*\omega_Y(aE) \rightarrow \mathcal{E}xt^{-n}(\shf{I}, \omega_{X}^{\bullet}),
\end{equation}
whose composition is the identity map.

By Claim (\ref{claim}) below, we have that $\mathcal{E}xt^{-n}(\shf{I}, \omega_{X}^{\bullet})\iso \omega_{X}$, thus (\ref{split of ext}) becomes
\begin{equation*}
\omega_{X}\rightarrow f_*\omega_Y(aE)\rightarrow \omega_{X},
\end{equation*}
so the natural morphism $f_*\omega_Y(aE)\rightarrow \omega_{X}$ is surjective and hence isomorphic.
\end{proof}

\begin{claim}\label{claim}
$\mathcal{E}xt^{-n}(\shf{I}, \omega_{X}^{\bullet})\iso \omega_{X}.$
\end{claim}
\begin{proof}
By definition of $\shf{I}$ we have the following sequence
\begin{equation*}
0\rightarrow \shf{I}\rightarrow \shf{O}_{X} \rightarrow \shf{O}_T \rightarrow 0,
\end{equation*}
where $T$ is a scheme with support equal to $f(E)$; in particular,
the codimension of $T$ is at least two. Applying $\mathcal{R}Hom_{X}( \cdot, \omega _{X}^{\bullet})$, we get
\begin{equation*}
\mathcal{R}Hom_{X}( \shf{O}_T, \omega _{X}^{\bullet})
\rightarrow \mathcal{R}Hom_{X}(\shf{O}_{X} , \omega _{X}^{\bullet})
\rightarrow \mathcal{R}Hom_{X}( \shf{I}, \omega _{X}^{\bullet})\xrightarrow{+},
\end{equation*}
that is
\begin{equation*}
\omega _T^{\bullet} \rightarrow \omega _{X}^{\bullet} \rightarrow \mathcal{R}Hom_{X}( \shf{I}, \omega _{X}^{\bullet})\xrightarrow{+}.
\end{equation*}
By taking $-n$th cohomology, we have the exact sequence
\begin{equation*}
h^{-n}(\omega _T^{\bullet})
\rightarrow \omega _{X}
\rightarrow \shf{E}xt^{-n}(\shf{I}, \omega _{X}^{\bullet})
\rightarrow h^{-n+1}(\omega _T^{\bullet}).
\end{equation*}
Note that the first and the last term of this exact sequence are zero
since $\dim T\le n-2$, which proves the claim.
\end{proof}

If a log resolution $f: Y\rightarrow X$ resolves the Jacobian ideal $\mathcal{J}ac_X$ of $X$, that is the extension $\mathcal{J}ac_X\cdot\sshf{Y}$ is invertible, then the smallest $a$ for which $f_*\omega_Y(aE)\iso \omega_X$ can be almost read off from Mather--Jacobian discrepancies, as discussed below.

We first recall some basics about Mather--Jacobian discrepancies and multiplier ideals, and refer to \cite{EIM2011} for details. Let $X$ be a variety (not necessarily normal) of dimension $n$ and $\Omega_X$ be the K\"ahler differential. Let $\nu: \widehat{X} \rightarrow X$
be the Nash blowup. As a closed subscheme of $\mathbb{P}(\bigwedge ^n \Omega_X)$ from its definition, $\widehat{X}$ has a line bundle
\begin{equation*}
\widehat{K}_X:=\shf{O}_{\mathbb{P}(\bigwedge ^n \Omega_X)}(1)| _{\widehat{X}}.
\end{equation*}
Any log resolution $f:Y \rightarrow X$ of the Jacobian ideal $\mathcal{J}ac_X$ factors through $\widehat{X}$
$$\xymatrix{
 Y\ar[r]_{\widehat{f}}\ar@/^1.0pc/[rr]^f & \widehat{X}\ar[r]_{\nu} &X}. $$
One thereby can define a divisor $D$ on $Y$ as
\begin{equation*}
D:= \widehat{K}_{Y/X}-J_{Y/X}=K_Y-\widehat{f}^*\widehat{K}_X-J_{Y/X},
\end{equation*}
where $\shf{O}_Y(-J_{Y/X})=\mathcal{J}ac_X\cdot \sshf{Y}$.

Put $D=\sum_i a_i(E_i, X)E_i$, where the summation is over the $f$-exceptional divisors. The coefficient $a_i(E_i, X)$ is called the Mather-Jacobian discrepancy of $E_i$. It depends only on the valuation of $K(X)$ determined by $E_i$, but not on the choice of $f$. 

Since $K_Y-\widehat{f}^*\widehat{K}$ is exceptional and effective, the direct image $f_*\sshf{Y}(D)$ is a sheaf of ideals on $X$, which is called the Mather-Jacobian multiplier ideal of $D$, and is denoted by $\widehat{\mathcal{J}}(X, \sshf{X})$. We will need the following local vanishing

\begin{theorem}({\cite[Theorem 3.5]{EIM2011}})\label{relative vanishing}
Given a sheaf of ideals $\frak{a}\subset \sshf{X}$. Let $f: Y\rightarrow X$ be a log resolution of $\mathcal{J}ac_X\cdot \frak{a}$. Let $Z_{Y/X}$ be an effective divisor such that $\frak{a}\cdot\sshf{Y}=\sshf{Y}(-Z_{Y/X})$. With the notations as above, for any $t\in \mathbb{R}_{\ge 0}$, we have
\begin{equation*}
    R^if_*\shf{O}_Y\left(\widehat{K}_{Y/X}-J_{Y/X}-\left\lfloor tZ_{Y/X}\right\rfloor\right)=0,
\end{equation*}
for all $i>0$.
\end{theorem}

\begin{theorem}\label{direct image}
Let $X$ be a normal variety and $f: Y\rightarrow X$ be a log resolution of $\mathcal{J}ac_X$. Then
$f_*\sshf{Y}\left(\widehat{f}^*\widehat{K}_X+J_{Y/X}\right)\iso\omega_X$.
\end{theorem}
\begin{proof}
Consider the following sequence
\begin{equation*}
0\rightarrow f_*\shf{O}_Y(D)\rightarrow \shf{O}_X \rightarrow i_*\shf{O}_W\rightarrow 0,
\end{equation*}
where $W$ is the subscheme of $X$ determined by the multiplier ideal $f_*\shf{O}_Y(D)$ and $i: W\hookrightarrow X$ is the inclusion. Since $W$ is supported on the singular locus of $X$ and $X$ is normal,
the codimension of $W$ is at least two.

Applying $\mathcal{R}Hom(\cdot, \omega _X ^{\bullet})$ to the above sequence we get the distinguished triangle
\begin{equation*}
\mathcal{R}Hom(i_*\shf{O}_W, \omega _X ^{\bullet})\rightarrow \mathcal{R}Hom(\shf{O}_X, \omega _X ^{\bullet})\rightarrow \mathcal{R}Hom(f_*\shf{O}_Y(D), \omega _X ^{\bullet})\xrightarrow{+1}.
\end{equation*}
Since $\mathcal{R}Hom(i_*\shf{O}_W, \omega _X ^{\bullet})\iso i_*\omega _W ^{\bullet}$, taking $-n$th cohomology yields the exact sequence
\begin{equation*}
h^{-n}(i_*\omega _W ^{\bullet})
\rightarrow \omega _X \rightarrow
\shf{E}xt^{-n}(f_*\shf{O}_Y(D), \omega _X ^{\bullet})
\rightarrow
h^{-n+1}(i_*\omega _W ^{\bullet})\rightarrow \dots,
\end{equation*}
where the first and the last cohomology sheaves are zero because $\dim W\le n-2$. So in particular, $\omega _X \iso
\shf{E}xt^{-n}(f_*\shf{O}_Y(D), \omega _X ^{\bullet})$.

To calculate this sheaf, we use Theorem \ref{relative vanishing} and Grothendieck duality
\begin{equation*}
\mathcal{R}Hom(f_*\shf{O}_Y(D), \omega _X ^{\bullet})\iso
\mathcal{R}Hom(\mathcal{R}f_*\shf{O}_Y(D), \omega _X ^{\bullet})\iso
\mathcal{R}f_*(\omega _Y^{\bullet}\otimes \shf{O}_Y(-D)).
\end{equation*}
Thus taking $-n$th cohomology gives the isomorphism
\begin{equation*}
\shf{E}xt^{-n}(f_*\shf{O}_Y(D), \omega _X ^{\bullet})\iso f_* (\omega _Y \otimes \shf{O}_Y(-D)).
\end{equation*}
Then the statement follows from 
\begin{equation*}
\omega _Y \otimes \shf{O}_Y(-D)\iso \shf{O}_Y\left(\widehat{f}^*\widehat{K}_X+J_{Y/X}\right).\qedhere
\end{equation*}
\end{proof}

\begin{corollary}\label{MJ Discrepancy}
Let $X$ be normal and $f:Y \rightarrow X$ be a log resolution of the Jacobian ideal $\mathcal{J}ac_X \subset \sshf{X}$.  Let $E$ denote the full exceptional divisor with components $E_i$ and $a(E_i, X)$ be the Mather-Jacobian discrepancy. If $a\ge \max_{i}\{-a(E_i, X)\}$, then it holds that
\begin{equation*}
    f_*\omega_Y(aE)\iso \omega_X.
\end{equation*}
\end{corollary}
\begin{proof}
Since $f$ factors through the Nash blowup, one can write
\begin{equation*}
K_Y=\widehat{f}^*\widehat{K}_X+J_{Y/X}+P-N,
\end{equation*}
where $P$ and $N$ are effective exceptional divisors without common components. By the assumption, we have $aE-N\ge 0$, therefore
\begin{equation*}
K_Y+aE\ge \widehat{f}^*\widehat{K}_X+J_{Y/X}.
\end{equation*}
Pushing it down, we deduce that $f_*\omega _Y(aE)\supseteq f_*\sshf{Y}(\widehat{f}^*\widehat{K}_X+J_{Y/X})\iso \omega _X$ by Theorem \ref{direct image}. This finishes the proof.
\end{proof}

\section{Depth Property}
In this section, we generalize Proposition \ref{Main} to arbitrary relative dimension.

\begin{theorem}\label{Main Theorem}
Let $f: Y\rightarrow X$ be a surjective proper morphism from a smooth variety $Y$ to a normal variety $X$. Set $d=\dim Y-\dim X$. Suppose $E=f^{-1}(X_{\text{sing}})$ is a divisor, with reduced induced closed subscheme structure. Then for any $m\ge 1$,
\begin{equation*}
    R^df_*\omega^{\otimes m}_{Y}(aE)/{\text{tor}}
\end{equation*}
is a reflexive sheaf if $a\gg 0$. 
\end{theorem}

As remarked in \cite[p.~171]{kollar86II}, the depth $R^if_*\omega_Y$ is one in general even if $X$ is smooth, so the analogous result for general $R^if_*$ does not hold for arbitrary $i$.

Keep the notation in \S 2.2, the result amounts to saying that the natural map $R^df_*\omega^{\otimes m}_{Y}(aE)\rightarrow j_*\paren{R^df_*\omega^{\otimes m}_V}$ is surjective when $a\gg 0$.

\begin{lemma}\label{$m$=1}
Theorem \ref{Main Theorem} holds for projective $X, Y$ and $m=1$.
\end{lemma}
\begin{proof}
Note to begin with that we can always assume that the considered $f$ has connected fibres by Stein factorization and Lemma \ref{depth under finite map}. 

We will do induction on the relative dimension $d$. The case $d=0$ is just Proposition \ref{Main}. Let $d>0$. Given $f: Y\rightarrow X$, fix a positive integer $a$ such that $R^df_*\omega_Y(bE)/\text{tor}$ stabilizes for $b\ge a$. Take a smooth, ample divisor $H$ on $Y$ such that $f|_H: H\rightarrow X$ is surjective, $H$ intersects $E$ transversally, and $R^if_*\omega_Y(aE+H)=0$ for all $i>0$ (the last condition can be made by using for instance \cite[Prop.~1.45, Prop.~2.69]{KollarMori08}). Consider the exact sequence
\begin{equation}\label{ses}
    \ses{\omega_Y}{\omega_Y(H)}{\omega_H}. 
\end{equation}

When $d=1$, we have
the surjection
\begin{equation*}
\pi_a:  f_*\omega_H(aE)\twoheadrightarrow R^1f_*\omega_Y(aE)/{\text{tor}}. 
\end{equation*}
For each $b>a$, the commutative diagram 
$$\xymatrix{
 f_*\omega_H(aE)\ar@{->>}[r]^{\pi_a}\ar@{^(->}[d]_{\varphi} & R^1f_*\omega_Y(aE)/{\text{tor}}\ar[d]^{\iso }\\
 f_*\omega_H(bE)\ar[r]^{\pi_b} &R^1f_*\omega_Y(bE)/{\text{tor}}, } $$
 implies that  the natural morphism 
 \begin{equation*}
\pi_b: f_*\omega_H(bE)\twoheadrightarrow R^1f_*\omega_Y(bE)/{\text{tor}}. 
\end{equation*}
surjects. 

Consider the  Stein factorization of $f|_H: H\xrightarrow{f'} X'\xrightarrow{\pi} X$, where $f'$ is birational, $X'$ is normal and $\pi$ is finite. If $b\gg 0$, it holds that 
\[f_*\omega_H(bE)\iso \pi_*f'_*\omega_H(bE)\iso \pi_*\omega_{X'}, \]
which is $S_2$ by Lemma \ref{depth under finite map}. Moreover, $\pi_*\omega_{X'}$ contains $\omega_X$ as a direct summand.  In fact, there exists a trace map $\text{tr}: \pi_*\sshf{X'}\rightarrow \sshf{X}$ since $X$ is normal. Applying $\shf{H}om_{\sshf{X}}(-, \omega_X)$ and the duality for finite map (eg. \cite[Prop. 5.67]{KollarMori08}), we get
$\omega_X\hookrightarrow \pi_*\omega_{X'}$, which gives rise to a splitting of the natural map $\pi_*\omega_{X'}\rightarrow \omega_X$. 

By sequence (\ref{filtration}) and the above consideration, for $b\gg 0$, we arrive at the diagram
$$\xymatrix{
 R^1f_*\omega_Y(bE)/{\text{tor}}\ar@{^(->}[r]\ar[rrd]& \omega_X\ar@{^(->}[r]& \pi_*\omega_{X'}\ar@{->>}[d]^{\pi_b} \\
 & & R^1f_*\omega_Y(bE)/{\text{tor}}, } $$
where the oblique composite map is indeed an isomorphism, because it is an endomorphism,  isomorphic over $U=X\backslash X_{\text{sing}}$. Therefore it follows that $R^1f_*\omega_Y(bE)/{\text{tor}}$ is a direct summand of $\pi_*\omega_{X'}$. Hence by Lemma \ref{S_k condition}, $R^1f_*\omega_Y(aE)/{\text{tor}}\iso R^1f_*\omega_Y(bE)/{\text{tor}}$ is $S_2$. 

Let $d>1$. Then (\ref{ses}) yields the isomorphism
\begin{equation}\label{iso of higher direct images}
    R^{d-1}f_*\omega_H(aE)\xrightarrow{\iso} R^df_*\omega_Y(aE).
\end{equation}
Again for each $b>a$, we have the commutative diagram
$$\xymatrix{
 R^{d-1}f_*\omega_H(aE)/{\text{tor}}\ar[r]^{\iso }\ar@{^(->}[d]_{\varphi} & R^df_*\omega_Y(aE)/{\text{tor}}\ar[d]^{\iso }\\
 R^{d-1}f_*\omega_H(bE)/{\text{tor}}\ar[r] &R^df_*\omega_Y(bE)/{\text{tor}}. } $$
Thus the left inclusion $\varphi$ gives rise to a splitting
\begin{equation*}
    R^{d-1}f_*\omega_H(bE)/{\text{tor}}\iso R^{d-1}f_*\omega_H(aE)/{\text{tor}}\oplus Q,
\end{equation*}
where $Q$ is the quotient of $\varphi$. By the induction hypothesis, for $b\gg 0$, $R^{d-1}f_*\omega_H(bE)/{\text{tor}}$ is $S_2$, so is $R^{d-1}f_*\omega_H(aE)/{\text{tor}}$ by Lemma \ref{S_k condition}. Combining this with (\ref{iso of higher direct images}), we finish the proof.
\end{proof}

\begin{proof}[Proof of Theorem \ref{Main Theorem}]
\vspace{0.5cm}

Step 1 (Reduction to projective case): Since the statement is local on $X$, by taking an affine cover, we can assume $X$ is quasi-projective. Then we can find a normal projective variety $\overline{X}$ with an open immersion $j: X\rightarrow \overline{X}$. By Nagata's compactification theorem \cite{Nagata63}, see also \cite{Conrad07}, there exists an open immersion $Y\rightarrow \overline{Y}$ of $\overline{X}$-schemes such that $\overline{Y}\rightarrow \overline{X}$ is proper. Let $\overline{E}$ be the Zariski closure of $E$ in $\overline{Y}$. Take an embedded log resolution $g: Z\rightarrow \overline{Y}$ such that $g$ is an isomorphism over $Y$ and the full $g$-exceptional divisor $E_g$ union $g^{-1}_*(\overline{E})$ is simple normal crossing. Thus we reach the commutative diagram
$$\xymatrix{
Y\ar[dd]_{f}\ar[r]^{j'}& Z\ar[d]^{g} \\
  & \overline{Y}\ar[d]^{\bar{f}}\\
 X\ar[r]^{j}&  \overline{X}.}$$
The diagram is indeed Cartesian by the properness of $f$ and the property of fibred product. Lastly note that $R^df_*\omega^{\otimes m}_Y(aE)\iso R^d(\bar{f}\circ g)_*\omega^{\otimes m}_Z(a(E_g+g^{-1}_*(\overline{E})))\big|_X$, thus we can assume $X$ is projective for our purpose.
\vspace{0.5cm}

Step 2: We can assume there exists an effective $f$-exceptional divisor $F$ on $Y$ such that $-F$ is $f$-ample. In fact, by a suitable sequence of blowups along smooth centers, one gets a proper birational morphism $\varphi: Y'\rightarrow Y$ such that there exists an effective $(f\circ\varphi)$-exceptional divisor $F$ on $Y'$ with the property that $-F$ is $(f\circ\varphi)$-ample, cf.\cite[Lemma 6.1]{HogadiXu09}. Let $E'$ be the full $(f\circ\varphi)$-exceptional divisor. Suppose $R^d(f\circ\varphi)_*\omega^{\otimes m}_{Y'}(bE')/{\text{tor}}$ is $S_2$ for some $b\ge 0$. Then if $a\gg 0$ such that
\begin{equation*}
    m\varphi^*K_Y+a\varphi^*E\ge mK_{Y'}+bE',
\end{equation*}
we have $R^d(f\circ\varphi)_*\omega^{\otimes m}_{Y'}(bE')/{\text{tor}}\subseteq R^df_*\omega^{\otimes m}_Y(aE)/{\text{tor}}$, so the latter sheaf is $S_2$.
\vspace{0.5cm}

Step 3 (cyclic covering): By properly choosing such $F$ in Step 2 and a sufficiently ample Cartier divisor $L$ on $X$, we can assume further that
\begin{enumerate}
  \item $-F+f^*L$ is ample \cite[Prop. 1.45]{KollarMori08},
  \item $\omega_Y\otimes \sshf{Y}(-F+f^*L)$ is globally generated.
\end{enumerate}
By Bertini theorem, for a general member $\Gamma\in |\omega^m_Y\otimes \sshf{Y}(-mF+mf^*L)|$, $\Gamma$ is smooth and $\Gamma$ intersects $E$ transversally.

Take an $m$-cyclic covering $\mu: \widehat{Y}\rightarrow Y$ with $\Gamma$ being the branch locus. We have that $\widehat{Y}$ is smooth and
\begin{equation*}
  \omega_{\widehat{Y}}(a\mu^*E) \iso \mu^*\sshf{Y}(mK_Y+aE+(m-1)(-F+f^*L)),
\end{equation*}
which implies that
\begin{equation*}
    R^d(f\circ\mu)_*\paren{\omega_{\widehat{Y}}(a\mu^*E)}\iso {\bigoplus}^m_{i=1} R^df_*\paren{\omega^{\otimes{i}}_Y(aE-(i-1)F)}\otimes\sshf{X}((i-1)L).
\end{equation*}
\vspace{0.3cm}

Step 4: Applying Lemma \ref{$m$=1} to $f_\circ \mu$, we deduce that $R^d(f_\circ \mu)_*\paren{\omega_{\widehat{Y}}(a\mu^*E)}/{\text{tor}}$ is $S_2$ when $a\gg 0$. 
Then the assertion follows from the above decomposition and Lemma \ref{S_k condition}.
%
\end{proof}


\bibliography{bibliography}{}
\bibliographystyle{alpha}

\end{document}